\documentclass[a4paper]{amsproc}

\usepackage{amssymb}
\usepackage{pstricks}
\usepackage{pstcol}
\usepackage{graphicx}
\usepackage{amsmath}
\usepackage{amsfonts}
\usepackage{latexsym}

\theoremstyle{plain}

\theoremstyle{definition}

\newtheorem{rem}{Remark}[section]

\numberwithin{equation}{section}

\title[Incomplete beta function]{The integrals in Gradshteyn and Ryzhik. \\
Part 11: The incomplete beta function}

\subjclass[2000]{Primary 33}


\author[K. Boyadzhiev]{Khristo N. Boyadzhiev}
\address{Department of Mathematics,
Ohio Northern University, Ada, OH 45810}
\email{k-boyadzhiev@onu.edu}

\author[L. Medina]{Luis A. Medina}
\address{Department of Mathematics,
Tulane University, New Orleans, LA 70118}
\email{lmedina@math.tulane.edu}

\author[V. Moll]{Victor H. Moll}
\address{Department of Mathematics,
Tulane University, New Orleans, LA 70118}
\email{vhm@math.tulane.edu}

\address{\hfill{\it Received 2 07
2008, revised ?? }\newline Departamento de Matem\'atica
\newline
Universidad T\'ecnica Federico Santa Mar\'{\i}a
\newline  Casilla 110-V,
\newline Valpara\'{\i}so, Chile}

\thanks{The third author wishes to acknowledge the partial support of  
NSF-DMS 0409968. The second author was partially supported as a graduate
student by the same grant.}

\begin{document}

{\begin{flushleft}\baselineskip9pt\scriptsize {\bf SCIENTIA}\newline
Series A: {\it Mathematical Sciences}, Vol. ?? (2009), ??
\newline Universidad T\'ecnica Federico Santa Mar{\'\i}a
\newline Valpara{\'\i}so, Chile
\newline ISSN 0716-8446
\newline {\copyright\space Universidad T\'ecnica Federico Santa
Mar{\'\i}a\space 2009}
\end{flushleft}}

\vspace{10mm} \setcounter{page}{1} \thispagestyle{empty}

\begin{abstract}
The table of Gradshteyn and Rhyzik contains some  integrals  that can be 
expressed in terms of the incomplete beta function. We describe some 
elementary properties of this function and use them to check some formulas 
in the mentioned table.
\end{abstract}

\maketitle

\newcommand{\nn}{\nonumber}
\newcommand{\ba}{\begin{eqnarray}}
\newcommand{\ea}{\end{eqnarray}}
\newcommand{\ift}{\int_{0}^{\infty}}
\newcommand{\ione}{\int_{0}^{1}}
\newcommand{\ifft}{\int_{- \infty}^{\infty}}
\newcommand{\no}{\noindent}
\newcommand{\realpart}{\mathop{\rm Re}\nolimits}
\newcommand{\imagpart}{\mathop{\rm Im}\nolimits}

\newtheorem{Definition}{\bf Definition}[section]
\newtheorem{Thm}[Definition]{\bf Theorem} 
\newtheorem{Example}[Definition]{\bf Example} 
\newtheorem{Lem}[Definition]{\bf Lemma} 
\newtheorem{Note}[Definition]{\bf Note} 
\newtheorem{Cor}[Definition]{\bf Corollary} 
\newtheorem{Prop}[Definition]{\bf Proposition} 
\newtheorem{Problem}[Definition]{\bf Problem} 
\numberwithin{equation}{section}

\maketitle

\section{Introduction} \label{intro} 
\setcounter{equation}{0}

The table of integrals \cite{gr} contains a large variety of definite 
integrals that involve the {\em incomplete beta } function  defined here by 
the integral 
\begin{equation}
\beta(a) = \int_{0}^{1} \frac{x^{a-1} \, dx }{1+x}.
\label{beta-def}
\end{equation}
\noindent
The convergence of the integral requires $a > 0$. Nielsen, who used this 
function extensively, attributed it to Stirling \cite{nielsen2}, page 17. 
The table \cite{gr} prefers 
to introduce first the {\em digamma function}
\begin{equation}
\psi(x) = \frac{d}{dx} \log \Gamma(x) = \frac{\Gamma'(x)}{\Gamma(x)},
\end{equation}
\noindent
and define $\beta(x)$ by the identity 
\begin{equation}
\beta(x) = \frac{1}{2} \left( \psi \left( \tfrac{x+1}{2} \right) - 
\psi \left( \tfrac{x}{2} \right) \right). 
\label{alt-def}
\end{equation}
\noindent
This definition appears
as $\mathbf{8.370}$ and (\ref{beta-def}) appears as $\mathbf{3.222.1}$. Here 
\begin{equation}
\Gamma(x) = \int_{0}^{\infty} t^{x-1} e^{-t} \, dt
\end{equation}
\noindent
is the classical gamma function. Naturally, both starting points for $\beta$
are equivalent, and Corollary \ref{coro-1} proves (\ref{alt-def}).  The 
value 
\begin{equation}
\gamma := -\psi(1) = -\Gamma'(1)
\end{equation}
\noindent
is the well-known {\em Euler's constant}. \\

In this paper we will prove elementary properties of this function and use 
them to evaluate some definite integrals in  \cite{gr}.

\section{Some elementary properties} \label{sec-elem} 
\setcounter{equation}{0}

The incomplete beta function admits a representation by series. 

\begin{Prop}
Let $a \in \mathbb{R}^{+}$. Then 
\begin{equation}
\beta(a) = \sum_{k=0}^{\infty} \frac{(-1)^{k}}{a+k}.
\end{equation}
\end{Prop}
\begin{proof}
The result follows from the expansion of $1/(1+x)$ in (\ref{beta-def}) 
as a geometric series. 
\end{proof}

\begin{Cor}
\label{coro-1}
The incomplete beta function is given by 
\begin{equation}
\beta(a) = \frac{1}{2} \left[ \psi \left( \frac{a+1}{2} \right) - 
\psi \left( \frac{a}{2} \right) \right].
\label{rela-11}
\end{equation}
\noindent
This is $\mathbf{8.370}$ in \cite{gr}.
\end{Cor}
\begin{proof}
The expansion for the digamma function $\psi$
\begin{equation}
\psi(t) = -\gamma - \sum_{k=0}^{\infty} \left( \frac{1}{t+k} - \frac{1}{k+1}
\right)
\end{equation}
\noindent
has been discussed in \cite{moll-gr10}. Then 
\begin{equation}
\psi \left( \frac{a}{2} \right)  = 
-\gamma - \sum_{k=0}^{\infty} \left( \frac{2}{a + 2k} - \frac{1}{k+1}
\right)
\end{equation}
\noindent
and 
\begin{equation}
\psi \left( \frac{a+1}{2} \right)  = 
-\gamma - \sum_{k=0}^{\infty} \left( \frac{2}{a + 2k+1} - \frac{1}{k+1}
\right). 
\end{equation}
\noindent
The identity (\ref{rela-11}) comes from adding these two expressions. 
\end{proof}

 These properties are now employed to prove some functional relations of the
incomplete beta function. The proofs will employ the identities 
\begin{eqnarray}
\psi(x+1) & = & \frac{1}{x} + \psi(x) \label{psi-1} \\
\psi(x) - \psi(1-x) & = & - \pi \cot( \pi x) \label{psi-2} \\
\psi(x + \tfrac{1}{2} ) - \psi(\tfrac{1}{2} -x) 
& = &  \pi \tan( \pi x) \label{psi-3} 
\end{eqnarray}
\noindent
that were established in \cite{moll-gr10}.  \\

\begin{rem}
Several of the evaluations presented here will employ the special values
\begin{equation}
\psi(n+1) = -\gamma + \sum_{k=1}^{n} \frac{1}{k},
\end{equation}
\noindent
that appears as $\mathbf{8.365.4}$, and 
\begin{equation}
\psi \left( \tfrac{1}{2} \pm n \right) = -\gamma + 
2 \left( \sum_{k=1}^{n} \frac{1}{2k-1} - \ln 2 \right),
\end{equation}
\noindent
that appears as $\mathbf{8.366.3}$.  \\

Many of the formulas in Section $\mathbf{4.271}$ employ the values 
\begin{equation}
\psi'(n) = \frac{\pi^{2}}{6} - \sum_{k=1}^{n-1} \frac{1}{k^{2}},
\label{formula-211}
\end{equation}
\noindent
that appear as $\mathbf{8.366.11}$ and  also $\mathbf{8.366.12/13}$:
\begin{equation}
\psi'( \tfrac{1}{2} \pm n ) = \frac{\pi^{2}}{2} \mp 
4 \sum_{k=1}^{n} \frac{1}{(2k-1)^{2}}. 
\label{formula-212}
\end{equation}

Higher order derivatives are given by 
\begin{eqnarray}
\psi^{(n)}(1) & = & (-1)^{n+1} n! \zeta(n+1)  \text{ and } \nonumber \\ 
\psi^{(n)}(\tfrac{1}{2}) & = & (-1)^{n+1} n! (2^{n+1}-1)\zeta(n+1). \nonumber
\end{eqnarray}
\end{rem}

\begin{Prop}
The incomplete beta function satisfies 
\begin{eqnarray}
\beta(x+1) & = & \frac{1}{x} - \beta(x), \label{beta-1}  \\
\beta(1-x) & = & \frac{\pi}{\sin \pi x} - \beta(x), \label{beta-2}  \\
\beta(x+1) & = & \frac{1}{x} - \frac{\pi}{\sin \pi x} + \beta(1-x).
\label{beta-3}
\end{eqnarray} 
\end{Prop}
\begin{proof}
Using (\ref{rela-11}) we have 
\begin{eqnarray}
\beta(x+1) & = & \frac{1}{2} \left[ \psi \left( \frac{x+2}{2} \right) - 
\psi \left( \frac{x+1}{2} \right) \right] = 
\frac{1}{2} \left[ \psi \left( \frac{x}{2} + 1 \right) - 
\psi \left( \frac{x+1}{2} \right) \right] \nonumber \\
        & = & \frac{1}{2} \left[ \frac{2}{x} + 
\psi \left( \frac{x}{2} \right) - 
 \psi \left( \frac{x+1}{2} \right) \right] \nonumber \\
 & = & \frac{1}{x} - \beta(x). \nonumber 
\end{eqnarray}
\noindent
This establishes (\ref{beta-1}).  To prove (\ref{beta-2}) we start with
\begin{eqnarray}
\beta(x) + \beta(1-x) & = & \frac{1}{2} \left[ 
\psi \left( \frac{1}{2} + \frac{x}{2} \right) - 
\psi \left( \frac{x}{2}  \right) + 
\psi \left( 1 -\frac{x}{2}  \right) - 
\psi \left( \frac{1}{2} - \frac{x}{2} \right)  \right]. 
\nonumber 
\end{eqnarray}
\noindent
The formula (\ref{beta-2}) now follows from (\ref{psi-2}) and (\ref{psi-3}).
\end{proof}

\section{Some elementary changes of variables} \label{sec-elemchan} 
\setcounter{equation}{0}

The class of integrals evaluated here are obtained from (\ref{beta-def}) 
by some elementary manipulations. 

\begin{Example}
\label{ex21}
The change $x = t^{p}$ in (\ref{beta-def}) yields
\begin{equation}
\beta(a) = p \int_{0}^{1} \frac{t^{ap-1} \, dt}{1+t^{p}}.
\end{equation}
\noindent
Replace $a$ by $\frac{a}{p}$ to obtain $\mathbf{3.241.1}$:
\begin{equation}
\int_{0}^{1} \frac{t^{a-1} \, dt}{1+t^{p}} 
= \frac{1}{p} \beta \left( \frac{a}{p} 
\right).
\label{32411}
\end{equation}
\end{Example}

\begin{Example}
The special case $p=2$ in Example \ref{ex21} gives 
\begin{equation}
\beta(a) = 2 \int_{0}^{1} \frac{t^{2a-1} \, dt}{1+t^{2}}.
\end{equation}
\noindent
Choose $a = \frac{b+1}{2}$, and relabel the variable of integration as $x$, 
to obtain $\mathbf{3.249.4}$:
\begin{equation}
\int_{0}^{1} \frac{ x^{b} \, dx}{1+x^{2}} =  \frac{1}{2} \beta \left( 
\frac{b+1}{2}  \right).
\end{equation}
\end{Example}

\begin{Example}
The evaluation of $\mathbf{3.251.7}$:
\begin{equation}
\int_{0}^{1} \frac{x^{a} \, dx}{(1+x^{2})^{2}} = -\frac{1}{4} + 
\frac{a-1}{4} \beta \left( \frac{a-1}{2} \right)
\label{32517}
\end{equation}
\noindent
comes from the change of variables $t = x^{2}$ and integration by 
parts. Indeed,
\begin{eqnarray}
\int_{0}^{1} \frac{x^{a} \, dx}{(1+x^{2})^{2}} & = & 
\frac{1}{2} \int_{0}^{1} t^{(a-1)/2} \, \frac{d}{dt} \frac{1}{1+t} \, dt 
\nonumber \\
& = & -\frac{1}{4} + \frac{a-1}{4} \int_{0}^{1} \frac{t^{(a-3)/2} \, dt}
{1+t} \, dt, \nonumber
\end{eqnarray}
\noindent
and (\ref{32517}) has been established. 
\end{Example}

\begin{Example}
Formula $\mathbf{3.231.2}$ states that 
\begin{equation}
\int_{0}^{1} \frac{x^{p-1} + x^{-p}}{1+x} \, dx = \frac{\pi}{\sin \pi p}.
\end{equation}
\noindent
The integrals is recognized as $\beta(p) + \beta(1-p)$ and its value follows
from (\ref{beta-2}). Similarly, $\mathbf{3.231.4}$ is
\begin{equation}
\int_{0}^{1} \frac{x^{p} - x^{-p}}{1+x} \, dx = \frac{1}{p} - \frac{\pi}
{\sin \pi p}. 
\end{equation}
\noindent
The integral is now recognized as $\beta(1+p) - \beta(1-p)$, and the 
result follows from (\ref{beta-3}).
\end{Example}

\begin{Example}
The evaluation of $\mathbf{3.244.1}$:
\begin{equation}
\int_{0}^{1} \frac{x^{p-1} + x^{q-p-1}}{1+x^{q}} \, dx = 
\frac{\pi}{q} \text{cosec} \frac{p \pi}{q}
\end{equation}
\noindent
is
\begin{equation}
I = \frac{1}{q} \left( \beta(p/q) + \beta(1 - p/q) \right)
\end{equation}
\noindent
according to (\ref{32411}). The result now follows from (\ref{beta-2}).
\end{Example}

\begin{Example}
The evaluation of $\mathbf{3.269.2}$:
\begin{equation}
\int_{0}^{1} x \, \frac{x^{p} - x^{-p}}{1+x^{2}} \, dx = \frac{1}{p} - \frac{\pi}{2 \sin( \pi p/2)}
\end{equation}
\noindent
is obtained by the change of variables $t = x^{2}$, that produces 
\begin{equation}
I = \frac{1}{2} \int_{0}^{1} \frac{t^{p/2} - t^{-p/2}}{1+t} \, dt =
\frac{1}{2} \left[ \beta \left( \frac{p}{2} + 1 \right) - 
 \beta \left( 1 -\frac{p}{2} \right) \right]. 
\end{equation}
\noindent
The result now follows from (\ref{beta-3}). 
\end{Example}

\section{Some exponential integrals} \label{sec-expo} 
\setcounter{equation}{0}

In this section we present some exponential integrals that may be 
evaluated in terms of the $\beta$-function. 

\begin{Example}
The  change of  variables $x = e^{-t}$ in (\ref{beta-def}) gives
\begin{equation}
\beta(a) = \int_{0}^{\infty} \frac{e^{-at} \, dt}{1+e^{-t}}. 
\label{expo-1}
\end{equation}
\noindent
This appears as $\mathbf{3.311.2}$ in \cite{gr}.
\end{Example}

\begin{Example}
The  evaluation of $\mathbf{3.311.13}$: 
\begin{equation}
\int_{0}^{\infty} \frac{e^{-px} + e^{-qx}}{1+e^{-(p+q)x}} \, dx = 
\frac{\pi}{p+q} \text{cosec}\left( \frac{\pi p }{p+q} \right)
\end{equation}
\noindent
is achieved by the change of  variables $t = (p+q)x$ that produces 
\begin{eqnarray}
I & = & \frac{1}{p+q} \int_{0}^{\infty} \frac{e^{-pt/(p+q)}}{1+e^{-t}} \, dt +
 \frac{1}{p+q} \int_{0}^{\infty} \frac{e^{-qt/(p+q)}}{1+e^{-t}} \, dt
\nonumber \\
& = & \frac{1}{p+q} \left[ \beta \left( \frac{p}{p+q} \right) + 
\beta \left( 1 - \frac{p}{p+q} \right) \right]. \nonumber
\end{eqnarray}
\noindent
The result now comes from (\ref{beta-3}).
\end{Example}

\section{Some trigonometrical integrals} \label{sec-trigo} 
\setcounter{equation}{0}

In this section we present the evaluation of some trigonometric integrals 
using the $\beta$-function. 

\begin{Example}
The change of variables $x = \tan^{2}t$ in (\ref{beta-def}) gives 
\begin{equation}
\beta(a) = 2 \int_{0}^{\pi/4} \tan^{2a-1}t \, dt.
\end{equation}
\noindent
Introduce the new parameter $b = 2a-1$ to obtain $\mathbf{3.622.2}$:
\begin{equation}
\int_{0}^{\pi/4} \tan^{b}t \, dt = \frac{1}{2} \beta \left( \frac{b+1}{2} 
\right).
\label{36222}
\end{equation}
\end{Example}

\begin{Example}
The  change of variables $x = \tan t$ in (\ref{32517}) gives
\begin{equation}
\int_{0}^{\pi/4} \tan^{a}t \, \cos^{2}t \, dt = 
-\frac{1}{4} + \frac{a-1}{4} \beta \left( \frac{a-1}{2} \right). 
\label{form-1}
\end{equation}
\noindent
Now use (\ref{beta-1}) to obtain
\begin{equation}
\beta \left( \frac{a-1}{2} \right) = \frac{2}{a-1} - 
\beta \left( \frac{a+1}{2} \right), 
\end{equation}
\noindent 
that converts (\ref{form-1}) to 
\begin{equation}
\int_{0}^{\pi/4} \tan^{a}t \, \cos^{2}t \, dt = 
\frac{1}{4} + \frac{1-a}{4} \beta \left( \frac{a+1}{2} \right).
\label{form-2}
\end{equation}
\noindent
This is the form in which $\mathbf{3.623.3}$ appears in \cite{gr}. Using this 
form and (\ref{36222}) we obtain $\mathbf{3.623.2}$:
\begin{equation}
\int_{0}^{\pi/4} \tan^{a}t \, \sin^{2}t \, dt = 
-\frac{1}{4} + \frac{1+a}{4} \beta \left( \frac{a+1}{2} \right).
\label{form-3}
\end{equation}
\end{Example}

\begin{Example}
The evaluation of $\mathbf{3.624.1}$:
\begin{equation}
\int_{0}^{\pi/4} \frac{\sin^{p}x \, dx}{\cos^{p+2}x} = \frac{1}{p+1}
\end{equation}
\noindent
can be done by writing the integral as
\begin{equation}
I = \int_{0}^{\pi/4} \tan^{p+2}x \, dx + \int_{0}^{\pi/4} \tan^{p}x \, dx.
\end{equation}
\noindent
These are evaluated using (\ref{36222}) to obtain
\begin{equation}
I = \frac{1}{2} \beta \left( \frac{p+3}{2} \right) + \frac{1}{2} 
 \beta \left( \frac{p+1}{2} \right). 
\end{equation}
\noindent
The rule (\ref{beta-1}) completes the proof. 
\end{Example}

\begin{Example}
The integral $\mathbf{3.651.2}$
\begin{equation}
\int_{0}^{\pi/4} \frac{\tan^{\mu}x \, dx}{1- \sin x \cos x} = 
\frac{1}{3} \left( \beta \left( \frac{\mu+2}{2} \right) + 
\beta \left( \frac{\mu+1}{2} \right) \right)
\label{36512}
\end{equation}
\noindent
can be established directly using the integral definition of $\beta$ given in 
(\ref{beta-def}). Simply observe that dividing the numerator and denominator 
of the integrand by $\cos^{2}x$ yields, after the change of 
variables $t = \tan x$, the identity
\begin{eqnarray}
\int_{0}^{\pi/4} \frac{\tan^{\mu}x \, dx}{1- \sin x \cos x}  & = & 
\int_{0}^{\pi/4} \frac{\tan^{\mu}x} {(\sec^{2}x- \tan x)}  \frac{dx}{\cos^{2}x} 
\nonumber  \\
& = & \int_{0}^{1} \frac{t^{\mu} \, dt}{t^{2}-t+1} \nonumber \\
& = & \int_{0}^{1} \frac{t^{\mu+1} + t^{\mu}}{t^{3}+1} \, dt. \nonumber 
\end{eqnarray}
\noindent
The change of variables $t = s^{1/3}$ gives the result.  \\

The evaluation of $\mathbf{3.651.1}$ 
\begin{equation}
\int_{0}^{\pi/4} \frac{\tan^{\mu}x \, dx}{1+ \sin x \cos x} = 
\frac{1}{3} \left( \psi \left( \frac{\mu+2}{2} \right) -
\psi \left( \frac{\mu+1}{2} \right) \right)
\label{36511}
\end{equation}
\noindent
can be established along the same lines.  This part employs the 
representation $\mathbf{8.361.7}$:
\begin{equation}
\psi(z) = \int_{0}^{1} \frac{x^{z-1}-1}{x-1} \, dx - \gamma
\end{equation}
\noindent
established in \cite{moll-gr10}.
\end{Example}

\begin{Example}
The elementary identity
\begin{equation}
\frac{1}{1 - \sin^{2}x \cos^{2}x} = \frac{1}{2} 
\left( \frac{1}{1+ \sin x \cos x} + \frac{1}{1 - \sin x \cos x} \right)
\end{equation}
\noindent 
and the evaluations given in Examples \ref{36511} and \ref{36512} gives a 
proof of $\mathbf{3.656.1}$: 
\begin{equation}
\tfrac{1}{12} \left( 
- \psi \left( \tfrac{\mu+1}{6} \right) 
- \psi \left( \tfrac{\mu+2}{6} \right) 
+ \psi \left( \tfrac{\mu+4}{6} \right) 
+ \psi \left( \tfrac{\mu+5}{6} \right) 
+ 2 \psi \left( \tfrac{\mu+2}{6} \right) 
- 2 \psi \left( \tfrac{\mu+1}{6} \right) \right).
\end{equation}
\end{Example}
	
\begin{Example}
The final integral in this section is $\mathbf{3.635.1}$:
\begin{equation}
\int_{0}^{\pi/4} \cos^{\mu-1}(2x) \, \tan x \, dx = \frac{1}{2} \beta(\mu).
\end{equation}
\noindent
This is easy: start with 
\begin{equation}
\tan x = \frac{\sin x}{\cos x} = \frac{2 \sin x \cos x}{2 \cos^{2}x} = 
\frac{\sin 2x}{1+ \cos 2x},
\end{equation}
\noindent
and use the change of variables $t = \cos 2x$ to produce the result. 
\end{Example}

\section{Some hyperbolic integrals} \label{sec-hyper} 
\setcounter{equation}{0}

This section contains the evaluation of some hyperbolic integrals using the
$\beta$-function. 

\begin{Example}
The integral (\ref{expo-1}) can be written as 
\begin{equation}
\beta(a) = \int_{0}^{\infty} \frac{e^{t(1/2-a)} \, dt}{e^{t/2} + e^{-t/2}},
\end{equation}
\noindent
and with $t = 2y$ and $b = 2a-1$, we obtain $\mathbf{3.541.6}$:
\begin{equation}
\int_{0}^{\infty} \frac{e^{-by} \, dy }{\cosh y} = \beta \left( \frac{b+1}{2} 
\right). 
\label{35416}
\end{equation}
\end{Example}

\begin{Example}
Integration by parts produces
\begin{eqnarray}
\int_{0}^{\infty} \frac{e^{-ax} \, dx}{\cosh^{2}x} & = & 
2 \int_{0}^{\infty} e^{-ax} \frac{d}{dx} \frac{1}{1+e^{-2x}} \, dx \nonumber \\
& = & -1 + 2a \int_{0}^{\infty} \frac{e^{-ax} \, dx}{1+ e^{-2x}}. \nonumber
\end{eqnarray}
\noindent
The change of variables $t = 2x$ now gives the evaluation of $\mathbf{3.541.8}$:
\begin{equation}
\int_{0}^{\infty} \frac{e^{-ax} \, dx}{\cosh^{2}x} = a \beta \left( 
\frac{a}{2} \right) -1. 
\end{equation}
\end{Example}

\begin{Example}
The change of variables $t = e^{-x}$ gives
\begin{equation}
\int_{0}^{\infty} e^{-ax} \tanh x \, dx = \int_{0}^{1} \frac{t^{a-1}-t^{a}}
{1+t^{2}} \, dt,
\end{equation}
\noindent
and with $s = t^{2}$ we get 
\begin{eqnarray}
I & = &  \frac{1}{2} \int_{0}^{1} \frac{s^{a/2-1} - s^{(a-1)/2}}{1+s} \, ds 
\nonumber \\
& = & \frac{1}{2} \left[ \beta \left( \frac{a}{2} \right) - 
\beta \left( \frac{a}{2} + 1 \right) \right].
\nonumber
\end{eqnarray}
\noindent
The transformation rule (\ref{beta-1}) 
gives the evaluation of $\mathbf{3.541.7}$: 
\begin{equation}
\int_{0}^{\infty} e^{-ax} \tanh x \, dx = \beta \left( \frac{a}{2} \right) 
- \frac{1}{a}.
\end{equation}
\end{Example}

\section{Differentiation formulas} \label{sec-diff} 
\setcounter{equation}{0}

\begin{Example}
Differentiating (\ref{beta-def}) with respect to the parameter $a$ yields
\begin{equation}
\int_{0}^{1} \frac{x^{a-1} \, \ln x}{1+x} \, dx = \beta'(a),
\end{equation}
\noindent
that appears as $\mathbf{4.251.3}$ in \cite{gr}. 
\end{Example}

\begin{Example}
Differentiating 
(\ref{32411}) $n$ times with respect to the parameter $a$ produces
$\mathbf{4.271.16}$ written in the form
\begin{equation}
\int_{0}^{1} \frac{x^{a-1} \, \ln^{n}x}{1+x^{p}} \, dx = 
\frac{1}{p^{n+1}} \beta^{(n)}\left( \frac{a}{p} \right). 
\label{427116}
\end{equation}
\noindent
The choice $n=1$ now gives formula $\mathbf{4.254.4}$ in \cite{gr}: 
\begin{equation}
\int_{0}^{1} \frac{x^{a-1} \, \ln x}{1+x^{p}} \, dx = 
\frac{1}{p^{2}} \beta'\left( \frac{a}{p} \right). 
\label{42544}
\end{equation}
\end{Example}

\begin{Example}
The special case $n=1, \, a=1$ and $p=1$ in (\ref{427116}) produces the 
elementary integral $\mathbf{4.231.1}$:
\begin{equation}
\int_{0}^{1} \frac{\ln x \, dx}{1+x} = - \frac{\pi^{2}}{12}.
\label{42311}
\end{equation}
\noindent
In this evaluation we have employed the values 
\begin{equation}
\psi'(1) = \zeta(2) = \frac{\pi^{2}}{6}, \text{ and }
\psi'(1/2) = \frac{\pi^{2}}{12},
\end{equation}
\noindent
that appear in (\ref{formula-212}).
\end{Example}

\begin{Example}
Formula $\mathbf{4.231.14}$:
\begin{equation}
\int_{0}^{1} \frac{x \, \ln x }{1+x^{2}} \, dx = - \frac{\pi^{2}}{48}
\end{equation}
\noindent
comes from (\ref{42544})
 by choosing the parameters $n=1, \, a=2$ and $p=2$. The values 
of $\psi'(1)$ and $\psi'(1/2)$ are employed again. Naturally, this evaluation
also comes from (\ref{42311}) via the change of variables $x^{2} \mapsto x$.
\end{Example}

\begin{Example}
The choice $n=a=1$ and $p=2$ in (\ref{42544})  and the values 
\begin{equation}
\psi^{(2)} \left( \tfrac{1}{4} \right) = \pi^{2} + 8G \text{ and }
\psi^{(2)} \left( \tfrac{3}{4} \right) = \pi^{2} - 8G, 
\end{equation}
\noindent
where $G$ is {\em Catalan constant} defined by 
\begin{equation}
G = \sum_{k=0}^{\infty} \frac{(-1)^k}{(2k+1)^{2}} 
\end{equation}
\noindent
yields the evaluation of $\mathbf{4.231.12}$:
\begin{equation}
\int_{0}^{1} \frac{\ln x \, dx}{1+x^{2}} = -G.
\end{equation}
\noindent
The change of variables $x = t/a$, with $a>0$, and the elementary integral
\begin{equation}
\int_{0}^{a} \frac{dt}{t^{2}+a^{2}} = \frac{\pi}{4a},
\end{equation}
\noindent
give the evaluation of $\mathbf{4.231.11}$:
\begin{equation}
\int_{0}^{a} \frac{\ln x \, dx}{x^{2}+a^{2}} = \frac{\pi \ln a - 4G}{4a}.
\end{equation}
\end{Example}

\begin{Example}
Now choose $n=1, \, a=2$ and $p=1$ in (\ref{42544}) and use the value
$\psi'(3/2) = \pi^{2}/2-4$ given in (\ref{formula-212}) 
to obtain $\mathbf{4.231.19}$:
\begin{equation}
\int_{0}^{1} \frac{x \, \ln x }{1+x} \, dx = \frac{\pi^{2}}{12} -1. 
\end{equation}
\noindent
Combining this with (\ref{42311}) gives $\mathbf{4.231.20}$:
\begin{equation}
\int_{0}^{1} \frac{1-x}{1+x} \, \ln x \, dx = 1 - \frac{\pi^{2}}{6}. 
\end{equation}
\end{Example}

\begin{Example}
The values 
\begin{equation}
\psi^{(2)} \left( \tfrac{1}{4} \right) = -2 \pi^{3} - 56 \zeta(3) \text{ and }
\psi^{(2)} \left( \tfrac{3}{4} \right) = 2 \pi^{3} - 56 \zeta(3),
\end{equation}
\noindent
given in \cite{srichoi}, are now used to produce 
the evaluation of $\mathbf{4.261.6}$:
\begin{equation}
\int_{0}^{1} \frac{\ln^{2}x \, dx}{1+x} = \frac{\pi^{3}}{16}.
\end{equation}
\end{Example}

\begin{Example}
The relation
\begin{equation}
\psi^{(n)}(1-z) + (-1)^{n+1} \psi^{(n)}(z) = (-1)^{n} \pi 
\frac{d^{n}}{dz^{n}} \cot \pi z,
\end{equation}
\noindent
and the choice $n=4, \, a=1$ and $p=2$ in (\ref{42544})  produces 
\begin{eqnarray}
\int_{0}^{1} \frac{\ln^{4}x \, dx}{1+x^{2}} & = & 
\frac{1}{2^{5}} \beta^{(4)} \left( \tfrac{1}{2} \right) \nonumber \\
& = & \frac{1}{1024} \left( \psi^{(4)} \left( \tfrac{3}{4} \right) - 
\psi^{(4)} \left( \tfrac{1}{4} \right) \right) \nonumber \\
& = & \frac{1}{1024} \left( - \pi \frac{d^{4}}{dz^{4}} \cot \pi z 
\Big|_{z=3/4} \right). \nonumber
\end{eqnarray}
\noindent
This yields the evaluation of 
$\mathbf{4.263.2}$:
\begin{equation}
\int_{0}^{1} \frac{\ln^{4}x \, dx }{1+x^{2}} = \frac{5 \pi^{5}}{64}. 
\end{equation}

The evaluation of $\mathbf{4.265}$:
\begin{equation}
\int_{0}^{1} \frac{\ln^{6}x \, dx }{1+x^{2}} = \frac{61 \pi^{7}}{256},
\end{equation}
\noindent
can be checked by the same method. 
\end{Example}

\begin{Example}
Now choose $n \in \mathbb{N}$ and take $a= n +1$ and 
$p=1$ in (\ref{42544}) to obtain the expression
\begin{equation}
I := \int_{0}^{1} \frac{x^{n} \, \ln^{2}x}{1+x} \, dx = 
\frac{1}{8} \beta^{(2)}(n+1). 
\end{equation}
\noindent
This is now expressed in terms of the $\psi-$function and then simplified
employing the relation 
\begin{equation}
\psi^{(m)}(z) = (-1)^{m+1} m! \zeta(m+1,z),
\end{equation}
\noindent
with the Hurwitz zeta function 
\begin{equation}
\zeta(s,z) := \sum_{k=0}^{\infty} \frac{1}{(z+k)^{s}}.
\end{equation}
\noindent
We conclude that
\begin{equation}
I = \frac{1}{4} \left( \zeta \left( 3, \frac{n+1}{2} \right) - 
                       \zeta \left( 3, \frac{n+2}{2} \right) \right). 
\end{equation}
\noindent
The elementary identity 
\begin{equation}
\zeta \left( s, \tfrac{a}{2} \right) - 
\zeta \left( s, \tfrac{a+1}{2} \right)  = 2^{s} \sum_{k=0}^{\infty} 
\frac{(-1)^{k}}{(k+a)^{s}},
\end{equation}
\noindent
is now used with $s=3$ and $a=n+1$ to obtain
\begin{equation}
\int_{0}^{1} \frac{x^{n} \, \ln^{2}x \, dx }{1+x} = 
2 \sum_{k=0}^{\infty} \frac{(-1)^{k}}{(k+n+1)^{3}}. 
\end{equation}
\noindent
This is finally transformed to the form 
\begin{equation}
\int_{0}^{1} \frac{x^{n} \, \ln^{2}x \, dx }{1+x} = 
(-1)^{n} \left( \frac{3}{2} \zeta(3) + 2 \sum_{k=1}^{n} 
\frac{(-1)^{k}}{k^{3}} \right).
\end{equation}
\noindent
This is $\mathbf{4.261.11}$ of \cite{gr}.   \\

The same method produces $\mathbf{4.262.4}$:
\begin{equation}
\int_{0}^{1} \frac{x^{n} \, \ln^{3}x \, dx }{1+x} = 
(-1)^{n+1} \left( \frac{7 \pi^{4}}{120} - 6 \sum_{k=0}^{n-1} 
\frac{(-1)^{k}}{(k+1)^{4}} \right).
\end{equation}
\end{Example}

\medskip

\begin{Example}
The method of the previous example 
yields the value of  $\mathbf{4.262.1}$:
\begin{equation}
\int_{0}^{1} \frac{\ln^{3}x \, dx}{1+x} = - \frac{7 \pi^{4}}{120}.
\end{equation}
\noindent
Here we use $\psi^{(3)}(1) = \pi^{4}/15$ and $\psi^{(3)}(1/2) = \pi^{4}$. \\

Similarly, $\psi^{(5)}(1) = 8 \pi^{6}/63$ and $\psi^{(5)}(1/2) = 8 \pi^{6}$
yields $\mathbf{4.264.1}$:
\begin{equation}
\int_{0}^{1} \frac{\ln^{5}x \, dx}{1+x} = - \frac{31 \pi^{6}}{252},
\end{equation}
\noindent
and $\psi^{(7)}(1) = 8 \pi^{8}/15$ and $\psi^{(7)}(1/2) = 136 \pi^{8}$
yields $\mathbf{4.266.1}$:
\begin{equation}
\int_{0}^{1} \frac{\ln^{7}x \, dx}{1+x} = - \frac{127 \pi^{8}}{240}.
\end{equation}
\end{Example}

\begin{Example}
A combination of the  evaluations given above produces $\mathbf{4.261.2}$:
\begin{equation}
\int_{0}^{1} \frac{\ln^{2}x \, dx}{1-x+x^{2}} = 
\frac{10 \pi^{3}}{81 \sqrt{3}}.
\end{equation}
\noindent
Indeed, 
\begin{eqnarray}
\int_{0}^{1} \frac{\ln^{2}x \, dx}{1-x+x^{2}} & = & 
\int_{0}^{1} \frac{1+x}{1+x^{3}} \ln^{2}x \, dx  \nonumber \\
& = & 
\int_{0}^{1} \frac{\ln^{2}x \, dx}{1+x^{3}} + 
\int_{0}^{1} \frac{x \, \ln^{2}x \, dx}{1+x^{3}}  \nonumber \\
& = & \frac{1}{27} \left( \beta^{(2)}(\tfrac{1}{3} ) + 
\beta^{(2)}(\tfrac{2}{3} )  \right) \nonumber \\
& = & \frac{1}{216} \left(  
\psi^{(2)} \left( \tfrac{2}{3} \right) - 
\psi^{(2)} \left( \tfrac{1}{3} \right) + 
\psi^{(2)} \left( \tfrac{5}{6} \right) -
\psi^{(2)} \left( \tfrac{1}{6} \right) \right) \nonumber \\
& = & \frac{\pi}{216} 
\left(\frac{d^{2}}{dz^{2}} \cot \pi z 
\Big|_{z=1/3} + 
\frac{d^{2}}{dz^{2}} \cot \pi z 
\Big|_{z=1/6} \right)  \nonumber \\
& = & \frac{\pi}{216} \left( \frac{8 \pi^{2}}{3 \sqrt{3}} + 8 \sqrt{3} \pi^{2} 
\right) = \frac{10 \pi^{3}}{81 \sqrt{3}}. \nonumber 
\end{eqnarray}
\end{Example}

\begin{Example}
Replace $n$ by $2n$ in (\ref{427116}) and set $a=p=1$ to 
produce 
\begin{eqnarray}
\int_{0}^{1} \frac{\ln^{2n}x \,  dx}{1+x} & = & \beta^{(2n)}(1) \nonumber \\
 & = & \frac{1}{2^{2n+1}} \left( \psi^{(2n)}(1) - \psi^{(2n)} \left( 
\tfrac{1}{2} \right) \right) \nonumber \\
& = & \frac{2^{2n}-1}{2^{2n}} (2n)! \zeta(2n+1). \nonumber 
\end{eqnarray}
\noindent 
This appears as $\mathbf{4.271.1}$. 
\end{Example}

\begin{Example}
The change of variables $t = bx$ in (\ref{427116}) produces 
\begin{eqnarray}
\int_{0}^{b} \frac{t^{a-1} \, \ln t }{b^{p} + t^{p}} & = & 
\frac{b^{a-p}}{p^{2}} \beta' \left( \tfrac{a}{p} \right) + 
b^{1-a} \ln b \int_{0}^{b} \frac{t^{a-1} \, dt}{b^{p} + t^{p}} \nonumber \\
 & = & \frac{b^{a-b}}{p^{2}} \beta' \left( \tfrac{a}{p} \right) + 
\ln b \frac{b^{a-p}}{p} \beta \left( \tfrac{a}{p} \right). \nonumber 
\end{eqnarray}
\noindent 
The last integral was evaluated using (\ref{32411}).  \\

Differentiate this identity with respect to the parameter $b$ to obtain
\begin{eqnarray}
\int_{0}^{b} \frac{t^{a-1} \, \ln t}{(b^{p}+t^{p})^{2}} \, dt & = & 
\frac{b^{a-2p} \ln b}{2p} + \frac{p-a}{p^{3}}b^{a-2p} 
\beta' \left( \tfrac{a}{p} \right) \label{mess-1} \\
&  - & \frac{b^{a-2p}}{p^{2}} ( 1 + (a-p) \ln b ) \beta \left( \tfrac{a}{p} 
\right). \nonumber 
\end{eqnarray}
\noindent
The special case $a=b=p=1$ yields $\mathbf{4.231.6}$:
\begin{equation}
\int_{0}^{1} \frac{\ln x \, dx}{(1+x)^{2}} = -\beta(1) = - \ln 2.  
\end{equation}
\noindent
Similarly, the choice $a=2, \, b=1$ and  $p=2$ yields $\mathbf{4.234.2}$: 
\begin{equation}
\int_{0}^{1} \frac{x \ln x \, dx}{(1+x)^{2}} = -\frac{1}{4} \beta(1) = 
- \frac{\ln 2}{4}.  
\end{equation}
\end{Example}

\begin{Example}
In this last example of this section 
we present an evaluation of $\mathbf{4.234.1}$:
\begin{equation}
\int_{1}^{\infty} \frac{\ln x \, dx}{(1+x^{2})^{2}} = 
\frac{G}{2} - \frac{\pi}{8},
\label{42341}
\end{equation}
\noindent
using the methods developed here.  We begin with the change of variables 
$x \mapsto 1/x$ to transform the problem to the interval $[0,1]$. We have 
\begin{equation}
\int_{1}^{\infty} \frac{\ln x \, dx}{(1+x^{2})^{2}} = 
- \int_{0}^{1} \frac{x^{2} \ln x \, dx}{(1+x^{2})^{2}}. 
\end{equation}
\noindent
Now choose $a=3, \, b=-1$ and $p=2$ in (\ref{mess-1}) to obtain 
\begin{equation}
\int_{0}^{1} \frac{x^{2} \ln x \, dx}{(1+x^{2})^{2}} = 
- \frac{1}{8} \beta' \left( \frac{3}{2} \right) + 
\frac{1}{4} \beta \left( \frac{3}{2} \right). 
\end{equation}
\noindent
The value of (\ref{42341})  now follows from 
\begin{eqnarray}
\tfrac{1}{4} \beta \left( \tfrac{3}{2} \right) & = & 
\tfrac{1}{8} \left( \psi \left( \tfrac{5}{4} \right) - 
\psi \left( \tfrac{3}{4} \right) \right) \nonumber \\
& = & \tfrac{1}{8} \left( 4 + \psi \left( \tfrac{1}{4} \right) - 
\psi \left( \tfrac{3}{4} \right) \right) \nonumber \\
& = & \tfrac{1}{2} - \tfrac{\pi}{8}, \nonumber 
\end{eqnarray}
\noindent
and  
\begin{eqnarray}
\tfrac{1}{8} \beta' \left( \tfrac{3}{2} \right) & = & 
\tfrac{1}{32} \left( \psi' \left( \tfrac{5}{4} \right) - 
\psi' \left( \tfrac{3}{4} \right) \right) \nonumber \\
& = & \tfrac{1}{32} \left( \psi' \left( \tfrac{1}{4} \right) - 
\psi' \left( \tfrac{3}{4} \right) - 16 \right) \nonumber \\
& = & \tfrac{1}{32} \left( \zeta(2, \tfrac{1}{4}) - 
\zeta(2, \tfrac{3}{4}) - 16 \right) \nonumber \\
& = & \tfrac{G}{2} - \tfrac{1}{2}. 
\nonumber
\end{eqnarray}
\end{Example}

\section{One last example} \label{sec-last} 
\setcounter{equation}{0}

In this section we discuss the evaluation of $\mathbf{3.522.4}$:
\begin{equation}
\ift \frac{dx}{(b^{2}+x^{2}) \cosh \pi x} = \frac{1}{b} \beta \left( b + 
\frac{1}{2} \right). 
\label{35224}
\end{equation}
\noindent
The technique illustrated here will be employed in a future publication to 
discuss many other evaluations. 

To establish (\ref{35224}), introduce the function
\begin{equation}
h(b,y) := \ift e^{-bt} \frac{\cos y t}{\cosh t} \, dt. 
\end{equation}
\noindent
This function is harmonic and bounded for $\realpart{b} >0$. Therefore it
admits a Poisson representation 
\begin{equation}
h(b,y) = \frac{1}{\pi} \int_{-\infty}^{\infty} h(0,u) 
\frac{b}{b^{2} +(y-u)^{2}} \, du. 
\end{equation}
\noindent
The value $h(0,u)$ is a well-known Fourier transform
\begin{equation}
h(0,u) = \ift \frac{\cos yt}{\cosh t} \, dt = 
\frac{\pi}{2 \, \cosh(\pi y/2)},
\end{equation}
\noindent
that appears as $\mathbf{3.981.3}$ in \cite{gr}. Therefore we have
\begin{equation}
h(b,y) = \frac{b}{2} \int_{-\infty}^{\infty} 
\frac{du}{\cosh(\pi u/2) \, \left[ b^{2} + (y-u)^{2} \right]}. 
\end{equation}
\noindent
The special value $y=0$ and (\ref{35416}) give the result (after replacing 
$b$ by $2b$ and $u$ by $2u$).  \\

\begin{Note}
Formula (\ref{35224}) can also be obtained by a direct contour integration. 
Details will be provided in a future publication.
\end{Note}

\medskip

We conclude with an interpretation of (\ref{35224}) in terms of the sine 
Fourier transform of a function related to $\beta(x)$. The proof is a simple
application of the elementary identity
\begin{equation}
\ift e^{xt} \sin bt \, dt = \frac{b}{b^{2}+x^{2}}.
\end{equation}
\noindent
The details are left to the reader. 

\begin{Thm}
Let 
\begin{equation}
\mu(x) := \ift \frac{e^{-xt} \, dt}{\cosh t} = 
\beta \left( \frac{x+1}{2} \right).
\end{equation}
\noindent
Then $\mathbf{3.522.4}$ in (\ref{35224}) is equivalent to the identity
\begin{equation}
\ift \mu(t) \sin xt \, dt = \mu \left( \frac{2x}{\pi} \right).
\end{equation}
\end{Thm}

\bibliography{../../../AllRef/biblio2}
\bibliographystyle{plain}
\end{document}